\documentclass[12pt]{amsart}
\usepackage{amssymb}
\usepackage{a4wide}
\usepackage[latin1]{inputenc}

\newcommand{\NN}{{\mathbb N}}
\newcommand{\Rad}[1]{{\mathrm{Rad}}\:(#1)}
\newcommand{\End}[1]{{\mathrm{End}}\:(#1)}
\newcommand{\Hom}[2]{{\mathrm{Hom}}\:(#1,#2)}
\newcommand{\Jac}[1]{{\mathrm{Jac}}\:(#1)}
\renewcommand{\Im}[1]{{\mathrm{Im}\:(#1)}}

\newtheorem{thm}{Theorem}
\newtheorem{cor}[thm]{Corollary}
\newtheorem{prop}[thm]{Proposition}

\begin{document}
\title{Rings whose modules are weakly supplemented are perfect.}
\date{\today}
\author{Engin Büyüka\c{s}ik}
\address{Izmir Institute of Technology, Izmir, Turkey}
\email{enginbuyukasik@iyte.edu.tr}

\author{Christian Lomp}
\address{Departamento de Matemática Pura da Faculdade de Ciências da Universidade do Porto, R.Campo Alegre 687, 4169-007 Porto, Portugal}
\email{clomp@fc.up.pt}

\thanks{
This paper was written while the first author was visiting the University of Porto. He wishes to thank the members of the Department of Mathematics for their kind hospitality and the Scientific and Technical Research Council of Turkey (T\"{U}BITAK) for their financial support. The second author was supported by {\it Funda{\c{c}\~ao} para a Ci\^encia e a Tecnologia} (FCT) through the {\it Centro de Matem\'atica da Universidade do
Porto } (CMUP)}

\begin{abstract}
In this note we show that a ring $R$ is left perfect if and only if every left $R$-module is weakly supplemented if and only if $R$ is semilocal and the radical of the countably infinite free left $R$-module has a weak supplement.
 \end{abstract}
 \maketitle

H.Bass characterized in \cite{Bass} those ring $R$ whose left $R$-modules have projective covers and termed them {\it left perfect rings}. He characterized them as those semilocal rings which have a left $t$-nilpotent Jacobson radical $\Jac{R}$. Bass' {\it semiperfect rings} are those whose finitely generated left (or right) $R$-modules have projective covers and can be characterized as those semilocal rings which have the property that idempotents lift modulo $\Jac{R}$. Kasch and Mares transfered in \cite{KaschMares} the notions of perfect and semiperfect rings to modules and characterized  semiperfect modules by a lattice theoretical condition as follows: a module $M$ is called {\it supplemented} if for any submodule $N$ of $M$ there exists a submodule $L$ of $M$ minimal with respect to $M=N+L$. The left perfect rings are then shown to be exactly those rings whose left $R$-modules are supplemented while the semiperfect rings are those whose finitely generated left $R$-modules are supplemented. Equivalently it is enough for a ring $R$ to be semiperfect if the left (or right) $R$-module $R$ is supplemented.
Recall that a submodule $N$ of a module $M$ is called {\it small}, denoted by $N\ll M$, if $N+L\neq M$ for all proper submodules $L$ of $M$. Weakening the ``supplemented''-condition one calls a module {\it weakly supplemented} if for every submodule $N$ of $M$ there exists a submodule $L$ of $M$ with $N+L=M$ and $N\cap L \ll M$. The semilocal rings $R$ are  precisely those rings whose finitely generated left (or right) $R$-modules are weakly supplemented. Again it is enough that $R$ is weakly supplemented as left (or right) $R$-module. Semilocal rings which are not semiperfect are examples of weakly supplemented modules which are not supplemented. In this note we prove that if $R$ is semilocal and the radical of the countably infinite free left $R$-module has a weak supplement, then $R$ has to be left perfect, i.e. every left $R$-module is supplemented.

Throught this note all rings are associative with unit and modules are considered to be unital. An ideal $I$ of a ring $R$ is called left $t$-nilpotent if for any family $\{a_i\}_{i\in\NN}$ of elements of $R$ there exists $n>0$ such that $a_1a_2\cdots a_n=0$. A ring $R$ is left perfect if and only if it is semilocal and $\Jac{R}$ is left $t$-nilpotent. Recall that an infinite family $\{A_\lambda \mid \lambda \in \Lambda\}$ of left ideals of $R$ is called {\it left vanishing} if given any sequence $a_1,a_2,\ldots, $ with $a_i \in A_{\lambda_i}$ and $\lambda_i\neq \lambda_j$ for all $i\neq j$, there exists a number $n\geq 1$ for which $a_1a_2a_3\cdots a_n=0$. Ware and Zelmanowitz proved in \cite[Theorem 1]{WareZelmanowitz} that for any endomorphism $f\in \End{F}$ of a free module $F$ and endomorphism which belongs to the Jacobson radical $\Jac{\End{F}}$, the family $\{ \pi_\lambda(\Im{f}) \}_\Lambda$ of left ideals of $R$ is left vanishing. Using this result we can prove our main theorem: 

\begin{thm}
The following statements are equivalent for a ring $R$:
\begin{enumerate}
\item[(a)] Every left $R$-module is weakly supplemented;
\item[(b)] $R^{(\NN)}$ is weakly supplemented;
\item[(c)] $R$ is semilocal and $\Rad{R^{(\NN)}}$ has a weak supplement in $R^{(\NN)}$.
\item[(d)] $R$ is left perfect.
\end{enumerate}
\end{thm}

\begin{proof} $(d)\Rightarrow (a) \Rightarrow (b) \Rightarrow (c)$ is clear and we just need to show $(c) \Rightarrow (d)$:
Set $F=R^{(\NN)}$ and denote $J=\Jac{R}$. Suppose that $R$ is semilocal and $JF=\Rad{F}$ has a weak supplement in $F$. Let $L$ be a weak supplement of $JF$ in $F$, i.e. $JF + L = F$ and $JF \cap L \ll F$. Then $R=\pi_i(JF+L) = J + \pi_i(L) = \pi_i(L)$ for any $i\in \NN$  implies that there exists $x_i \in L$ such that $\pi_i(x_i)=1$.
Let $\{ a_i \}_{i\in \NN}$ be any family of elements of $J$ then $a_ix_i \in JL \subseteq JF \cap L \ll F$ and $\pi_i(a_ix_i)=a_i$ for any $i\in \NN$. Define $f:F\rightarrow F$ by $f(z)=\sum_{i\in \NN} z_ia_ix_i$ for all $z\in F$. Since $\Im{f} \ll F$, we get by Ware and Zelmanowitz's Theorem \cite[Theorem 1]{WareZelmanowitz} that
$\{ \pi_i(JL)\}_{i\in\NN}$ is left vanishing. Thus there exists $n>0$ such that 
$$a_1a_2\cdots a_n = \pi_1(a_1x_1)\pi_2(a_2x_2)\cdots \pi_n(a_nx_n)=0.$$
This shows that $\Jac{R}$ is left $t$-nilpotent and hence $R$ is left perfect.
\end{proof}

Let $\sigma[M]$ denote the Wisbauer category of a module $M$, i.e. the full category of $R$-Mod consisting of submodules of quotients of direct sums of copies of $M$. A module $M$ is called a self-generator if any of its  submodules is an image of a direct sum of copies of $M$.

\begin{cor}
Let $M$ be a finitely generated, self-projective, self-generator. Then every module in $\sigma[M]$ is weakly supplemented if and only if $\End{M}$ is left perfect.
\end{cor}

\begin{proof}
 By \cite[18.3]{Wisbauer} $M$ is projective in $\sigma[M]$ and by \cite[8.5]{Wisbauer} $M$ is a generator in $\sigma[M]$. Hence by \cite[46.2]{Wisbauer} the functor $\Hom{M}{-}$ is a Morita equivalence between $\sigma[M]$ and $\End{M}$-Mod. Thus every module in $\sigma[M]$ is  weakly supplemented if and only if every left $\End{M}$-module is weakly supplemented, which holds if and only if $\End{M}$ is left perfect by the Theorem.
\end{proof}

We finish the paper with a comment on weak supplements of images of endomorphisms. Recall that a left $R$-module $M$ is called semi-projective if for any endomorphism $f\in S=\End{M}$ we have $Sf = \Hom{M}{\Im{f}}$. The module $M$ is called $\pi$-projective if for any submodules $N,L$ of $M$ with $M=N+L$ we have $S=\Hom{M}{N}+\Hom{M}{L}$.

\begin{prop}
Suppose $M$ is a semi-projective and $\pi$-projective $R$-module. Then $S/\Jac{S}$ is regular if and only if $\Im{f}$ has a weak supplement in $M$ for each $f \in S$.
\end{prop}

\begin{proof} $(\Rightarrow)$ Let $f \in S$. By hypothesis there is
a $g \in S$ such that $f-fgf \in J(S)$. We have $\Im{f} + \Im{1-fg}=M$. It is easy to see that $\Im{f} \cap \Im{1-fg}\subseteq \Im{f-fgf}$, but since $f-fgf \in \Jac{S}$ we have $\Im{f-fgf} \ll M$. Hence $\Im{1-fg}$ is a weak supplement of $\Im{g}$ in $M$. 

$(\Leftarrow )$ Let $f \in S$ and $K$ be a weak supplement of $\Im{f}$ in $M$. Since $M$ is semi-projective and $\pi$-projective we have $S=\Hom{M}{\Im{f}} + \Hom{M}{K} = Sf + \Hom{M}{K}$. Since $Sf \cap \Hom{M}{K} = \Hom{M}{\Im{f}\cap K}$ and $\Im{f}\cap K \ll M$, we get $Sf \cap \Hom{M}{K}\subseteq \Jac{S}$. Thus $Sf$ has weak supplement for all $f$, which implies $S/\Jac{S}$ being von Neumann regular by \cite[3.18]{Lomp}.
\end{proof}

The last proposition generalizes \cite[3.18]{Lomp}. Also as a consequence we conclude that the endomorphism ring of a semi-projective, $\pi$-projective weakly supplemented module is regular modulo its Jacobson radical.

\end{document}